\documentclass[12pt]{amsart}

\usepackage{amsthm}
\usepackage{amsfonts}
\usepackage{graphicx}

\newcommand{\R}{{\mathbb{R}}}

\newcommand{\N}{{\mathbb{N}}}
\newcommand{\D}{{\bf{D}}}

\newcommand{\E}{{{\bf E}}}
\newcommand{\Af}{{\mathcal{A}}}
\newcommand{\Bf}{{\mathcal{B}}}
\newcommand{\Cf}{{\mathcal{C}}}
\newcommand{\Df}{{\mathcal{D}}}
\newcommand{\Ef}{{\mathcal{E}}}

\newcommand{\Hf}{{\mathcal{H}}}
\newcommand{\Sf}{{\mathcal{S}}}
\newcommand{\Tf}{{\mathcal{T}}}

\newcommand{\Wf}{{\mathcal{W}}}
\newcommand{\talpha}{{\tilde{\alpha}}}
\newcommand{\tbeta}{{\tilde{\beta}}}
\newcommand{\tomega}{{\tilde{\omega}}}

\newtheorem{theorem}{Theorem}[section]

\newtheorem{proposition}[theorem]{Proposition}
\newtheorem{conjecture}[theorem]{Conjecture}

\newtheorem{definition}[theorem]{Definition}

\begin{document}
\title{On the Heegaard genus of $3$-manifolds obtained by gluing}
\author{Trent Schirmer}

\begin{abstract}
We introduce a new technique for finding lower bounds on the Heegaard genus of a $3$-manifold obtained by gluing a pair of $3$-manifolds together along an incompressible torus or annulus.  We deduce a number of inequalities, including one which implies that $t(K_1\# K_2)\geq \max \{t(K_1),t(K_2)\}$, where $t(-)$ denotes tunnel number, $K_1$ and $K_2$ are knots in $S^3$, and $K_1$ is $m$-small.  This inequality is best possible.  We also provide an interesting collection of examples, similar to a set of examples found by Schultens and Wiedmann \cite{schul-weid}, which show that Heegaard genus can stay persistently low under the kinds of gluings we study here.
\end{abstract}

\maketitle

\section{Introduction}

In this paper we study the way that Heegaard genus behaves when one glues two manifolds $M_1$ and $M_2$ together along annuli or tori in their boundary to obtain a new manifold $M$.  In particular, we ask how low the Heegaard genus of $M$ can be compared to the Heegaard genera of $M_1$ and $M_2$.  This is a questions which is easy to ask but difficult to answer, and it has received a great deal of attention.\\

It is a well known result of Haken's \cite{Hak} that every Heegaard surface in a connected sum $M=M_1\# M_2$ (which can be regarded as a gluing along spheres) arises naturally as the connected sum of Heegaard surfaces $F_1$ and $F_2$ for $M_1$ and $M_2$, respectively, and thus $g(M)=g(M_1)+g(M_2)$, where $g(-)$ denotes Heegaard genus.  The same can be said for boundary connect sums $M_1\#_\partial M_2$, which is an attachment of $M_1$ to $M_2$ along disks in their boundary.  However, as soon as we make the leap to gluings along surfaces of non-positive Euler characteristic, the situation becomes more complicated and the volume of literature on the subject is large.  We briefly touch on a few developments that are of relevance to our work here, without making any claim to completeness.\\

In the case that $M_1$ and $M_2$ are glued together along closed genus $n$ surfaces in their boundary, it is easy to see that $g(M)\leq g(M_1)+g(M_2)-n$ because one can simply amalgamate any generalized Heegaard splitting induced by Heegaard splittings of $M_1$ and $M_2$.  A great deal of results have been published which give conditions on the gluing map that ensure that such splittings are the only minimal genus splittings of $M$, so that this upper bound is actually attained. Scharlemann and Tomova's paper \cite{schartom} is a good early example of the combinatorial approach, Lackenby's elegant paper \cite{Lack} is a good example of the geometric approach to this problem, and there are many, many more recent papers on the subject which are too numerous to name here. \\

On the other hand, it is often the case that Heegaard genus will be much lower.  For example, if we let $M_1$ be a handlebody, and $M_2$ be the closure of the complement of a regular neighborhood of a graph in $S^3$, we see that it is also possible to have a situation where $g(M)=g(M_2)-g(F)=0$.  Lower bounds on $g(M)$ have been found, however.  Important examples include those found by Schultens \cite{schul2}, Lei and Yang \cite{chinese}, and Derby-Talbot\cite{talbot}, and the latter's main theorem has hypotheses which correspond closely to the hypotheses of the main theorem here.   As far as existence is concerned, in the more restrictive case that the gluing occurs along a pair of {\em incompressible} tori, Schultens and Wiedmann have found examples where $g(M)<(g(M_1)+g(M_2))/2$, and more specifically $g(M)=\max \{g(M_1),g(M_2)\}-1$.\\

The extensive literature on the behavior of tunnel number of knots in $S^3$ under connected sum also falls under our umbrella here (Moriah's survey \cite{ymor} is a good place to start).  The tunnel number of knot $K$ can be defined by $t(K)+1=g(E(K))$, where $E(K)$ is the exterior of $K$, and $E(K_1\# K_2)$ is the annular sum of $E(K_1)$ and $E(K_2)$ along meridional annuli.  Here it is not difficult to see that $t(K_1\# K_2)\leq t(K_1)+t(K_2)+1$, and analogously to the case of gluings along closed surfaces, a great deal of work has gone into showing that this bound can be acheived and the circumstances under which it is.  Here papers by Moriah and Rubinstein \cite{morrub}, Morimoto Sakuma and Yokota \cite{mor2}, and Kobayashi and Rieck \cite{kobri} are a few of the earlier papers establishing existence, and each uses a distinct approach.  The work of Kobayashi and Rieck in particular plays an important role in the construction of our links in Section 6 below.\\

On the lower end, the first example of a pair of knots in $S^3$ for which $t(K_1\# K_2)< t(K_1)+t(K_2)$ was found by Morimoto \cite{mor3}, and Kobayashi \cite{kob} soon after found that the difference $t(K_1)+t(K_2)-t(K_1\# K_2)$ can be made arbitrarily large.  Lower bounds on subadditivity were found by Schultens and Scharlemann \cite{schulschar}\cite{schulschar2}, and \cite{schulschar} in particular is of interest to us here, because there lower bounds are established on the {\em degeneration ratio} of two (or more) knots, defined by the formula $$d(K_1,K_2)=\frac{t(K_1)+t(K_2)-t(K_1\# K_2)}{t(K_1)+t(K_2)}.$$  Scharlemann and Schultens found that $d(K_1,K_2)\leq \frac{3}{5}$ for any prime knots $K_1$ and $K_2$ in $S^3$, and in the case that $E(K_1\# K_2)$ admits a strongly irreducible minimal genus Heegaard splitting, $d(K_1,K_2)\leq \frac{1}{2}$.\\

In the following paper, our main goal is to prove lower bounds on subadditivity of Heegaard genus under toroidal and annular gluings.  Section 2 is devoted to setting down notation, definitions, and foundational propositions which the expert will find familiar.  Section 3 is devoted to the development of what might be called substitution machinery, and plays an essential in the proof of Section 4.\\

The main result of this paper is proved in Section 4.  It shows that if $M$ is obtained by gluing $M_1$ to $M_2$ along a pair of incompressible tori $T_1\subset \partial M_1$, $T_2\subset \partial M_2$, and if $M$ admits a minimal genus Heegaard surface $F$ which nicely (in the sense of satisfying Properties A and B of Section 4) intersects the incompressible torus $T\subset M$ which is the image of $T_1$ and $T_2$ after gluing, then $g(M)\geq \max \{g(M_1),g(M_2)\}-1$.  The examples of Schultens and Wiedmann show this bound to be best possible.    The philosophical approach the proof is simple, although its execution is a bit technical.  The idea is that, under sufficiently nice circumstances, the surface $F\cap M_2$ will cut $M_2$ into compression bodies, and in a nice pattern.  We then take a close look at this pattern and create a simpler {\em doppelg\"{a}nger} surface $Q$ properly embedded in the solid torus $W$ which cuts it into handlebodies and also has a pattern corresponding to the one left by $F\cap M_2$ in $M_2$.  We then attach $W$ to $M_1$ instead of $M_2$, and in such a way that $Q$ matches up nicely with $F\cap M_1$ after the gluing to form a closed Heegaard surface $F'$ in $M_1\cup W$ of genus less than or equal to $F$.  After removing a core of $W$ we get back $M_1$, and after at most one stabilization of $F'$, we obtain a Heegaard splitting of the desired kind.\\

In Section 5 we give applications of the main theorem.  In the case when $M$ is obtained by a toroidal gluing of $M_1$ and $M_2$ along what we have called an $\Ef$-incongruous gluing map, we show that $g(M)\geq \min\{g(M_1),g(M_2)\}-1$.  The $\min$ term arises, as opposed to the $\max$ of the general theorem, due to a lack of control over precisely which slope a Heegaard surface might intersect the gluing torus.  The main result is then applied with better force to the case of a connected sum between knots, which can be interpreted as a very special form of toroidal gluing that allows us to control the slope in which Heegaard surfaces intersect the gluing torus.  Here we deduce that $t(K_1\# K_2)\geq \max\{t(K_1),t(K_2)\}$ when either of $K_1$ or $K_2$ is $m$-small, or whenever a minimal genus Heegaard splitting of $E(K_2\# K_2)$ can be made to intersect the decomposing annulus of the connected sum essentially.  The latter assumption is always realized when $E(K_1\# K_2)$ admits a strongly irreducible minimal genus splitting.\\

Section 6 is then devoted to small observations and conjectures which we believe likely to be attainable using methods similar to those employed here. It also contains a construction of an interesting class of examples of manifolds whose Heegaard genus persistently decreases under toroidal gluing.  Specifically for all $n>0$ we construct $n$-component links $L=L_1\cup \cdots \cup L_n\subset S^3$, $n$ knots $K_1,\cdots K_n$ in $S^3$, and homeomorphisms $h_i:\partial N(L_i)\rightarrow \partial N(K_i)$ such that $g(M)=g(E(L))-n$, where $M=E(L)\cup_{\cup h_i}(\bigcup E(K_i))$.  And in fact it is easy to see that if the gluing is performed in sequence, so that $M_0=E(L)$, $M_1=E(L)\cup_{h_1}E(K_1)$, and more generally $M_k=E(L)\cup_{h_1\cup \cdots \cup h_{k}}(E(K_1)\cup \cdots \cup E(K_{k}))$, then $g(M_i\cup_{h_i} E(K_{i}))=g(M_i)-1=\max\{g(M_i),g(E(K_i))\}-1$.  This does not improve the result of Schultens and Wiedmann.  It is distinct from their examples, which use Seifert fibered spaces where we use non-trivial, non-torus link complements in $S^3$, although it is similar because both constructions utilize bridge spheres of the knots $K_i$ in the same fashion.\\

\begin{center}
{\bf Acknowledgments}
\end{center}

I would like to acknowledge Charlie Frohman, my advisor Maggy Tomova, and Jesse Johnson for many helpful and encouraging conversations.

\section{Preliminaries}

We assume an elementary level of knowledge regarding $3$-manifold topology and knot theory as presented in \cite{Hem},\cite{Jac}, and \cite{Rolfsen}.  Moreover, the more specialized material on Heegaard splittings that occurs here is only very tersely presented; we recommend \cite{Schar1} and \cite{saitoschulschar} to unfamiliar readers.\\  

We work in the PL category and we assume our manifolds to be oriented and our gluing maps to be orientation preserving, although we rarely refer to orientations unless we need to, or when special emphasis is warranted.  A closed regular neighborhood of a polyhedron $Y$ embedded in a manifold $X$ is denoted $N(Y,X)$, and unless specified otherwise we implicitly assume that choices are made to be small and mutually consistent with one another, so that for any other polyhedron $Z\subset X$ occuring in our argument, $N(Y,X)\cap Z=N(Y\cap Z, Z)$ (particularly, if $Y\cap Z=\emptyset$, $N(Y,X)\cap Z=\emptyset$).  When $X$ is understood we will write simply $N(Y)$ for $N(Y,X)$.  With analogous conventions, let $E(Y,X)=\overline{X\setminus N(Y,X)}$ denote the {\em exterior} of $Y$ in $X$, and define the {\em frontier} $Fr(Y)$ of $Y$ to be $N(Y,X)\cap E(Y,X)$.  Unless otherwise specified we assume that embedded polyhedra meet one another in general position and that all embeddings of manifolds are proper, meaning that if $Y$ is a manifold embedded in $X$, $\partial Y$ is embedded in $\partial X$.\\  

If $X$ is a topological space with connected components $X_1,\cdots , X_n$, then $\hat{X}$ will denote the set $\{X_1,\cdots , X_n\}$, and $|X|$ will denote the cardinality of $\hat{X}$. If two components of $E(Y,X)$ meet a common component $N(Y')$ of $N(Y)$, we say they are {\em adjacent along $Y_i$}. Finally, if $Y$ is a separating $n$ manifold properly embedded in an $n+1$-manifold $X$, the {\em connectivity graph} $\Cf_Y$ is the graph which has one vertex corresponding to each component of $E(Y,X)$ and an edge corresponding to each component $Y'$ of $Y$ connecting the pair of vertices whose corresponding components are adjacent along $Y'$.\\

We now introduce a selection of definitions and some important propositions pertaining to Heegaard splittings which are well-known to experts.\\

\begin{definition}
A {\em compression body} $V$ is either a ball or is obtained from a thickened closed surface $F\times I$ of positive genus by attaching $2$-handles to $F\times \{0\}$ and $3$-handles along any resulting spherical boundary components.  We define $\partial_+V=F\times\{1\}$ and $\partial_-V=\partial V \setminus \partial_+V$.\\
\end{definition}

\begin{definition}
A {\em spine} of a compression body $V$ is a graph $X$ properly embedded in $V$ so that its univalent vertices lie in $\partial_-V$, and so that $E(X,V)\cong \partial_+V\times I$.  A {\em core } of $V$ is a disjoint union of simple closed curve embedded in a spine of $V$.\\
\end{definition}

\begin{definition}
A disjoint collection $\D$ of pairwise non-isotopic compressing disks for a compression body $V$ is {\em complete} if each component of $E(\D,V)$ is a ball or thickened surface.\\
\end{definition}

\begin{proposition}
Let $P$ be an incompressible surface properly embedded in a compression body $V$.  Then there exists a complete collection $\D$ of compressing disks for $V$ which intersects $P$ only in essential arcs, if at all.\\

\end{proposition}

\begin{proposition}
If $V$ is a compression body and $P$ is an incompressible surface properly embedded in $V$ with $\partial P\subset \partial_+V$, then $E(P,V)$ is a union of compression bodies.\\

\end{proposition}

\begin{definition}
An annulus $A$ properly embedded in a compression body $V$ is {\em spanning} if it has one component on $\partial_+V$ and the other on $\partial_-V$.\\
\end{definition}

\begin{proposition}
If $V$ is a compression body and $P$ is a connected, incompressible, $\partial$-incompressible surface properly embedded in $V$, then $P$ is either a disk, a spanning annulus, or isotopic into $\partial_-V$.\\
\end{proposition}

\begin{definition}

A {\em Heegaard splitting} $(V,W,F)$ of a manifold $M$ is a decomposition $M=V\cup W$ where each of $V,W$ is a compression body and $F=\partial_+ V=\partial_+ W=V\cap W$.  $F$ is called a {\em Heegaard surface}.  The {\em Heegaard genus} $g(M)$ of a manifold is the minimal genus of a Heegaard surface for $M$.\\
\end{definition}

\begin{definition}
A Heegaard splitting $(V,W,F)$ is said to be:

\begin{itemize}
\item {\em Stabilized} if there exist compressing disks $D\subset V$, $D'\subset W$ such that $D\cap D'$ consists of exactly one point.
\item {\em Reducible} if there exist compressing disks $D\subset V$, $D'\subset W$ such that $\partial D =\partial D'$.
\item {\em Weakly reducible} if there are compressing disks $D\subset V$, $D'\subset W$ such that $D\cap D' =\emptyset$
\item {\em Strongly irreducible} if, for all compressing disks $D\subset V$, $D'\subset W$, $D\cap D' \neq \emptyset$
\end{itemize}
\end{definition}

Every stabilized Heegaard splitting is reducible, and every reducible Heegaard splitting is weakly reducible.  By Haken's lemma, only irreducible manifolds admit weakly reducible or strongly irreducible Heegaard splittings, and it is a result of Casson and Gordon \cite{Cas-Gor} that every manifold which admits a weakly reducible Heegaard splitting is Haken. Subsequently Scharlemann and Thompson \cite{scharthomp} developed the theory of thin position in generalized Heegaard splittings.

\begin{definition}
A {\em generalized Heegaard splitting} $((V_1,W_1,F_1),\cdots , (V_n,W_n,F_n))$ of a manifold $M$ is a decomposition $M=M_1\cup\cdots \cup M_n$ such that $(V_i,W_i,F_i)$ forms a Heegaard splitting for the submanifold $M_i$, $M_i\cap M_{i+1}=\partial_-W_i=\partial_-V_{i+1}=S_i$ for $1\leq i<n$, and $M_i\cap M_j=\emptyset$ whenever else $i\neq j$.  The surfaces $F_i$ are called the {\em thick surfaces} of the generalized Heegaard splitting, while the surfaces $S_i$ are called {\em thin surfaces}.\\

\end{definition}

There is process of {\em untelescoping} whereby a weakly irreducible Heegaard splitting $(V,W,F)$ is changed into a generalized Heegaard splitting $((V_1,W_1,F_1),\cdots , (V_n,W_n,F_n))$ satisfying $g(F)=\Sigma g(F_i) -\Sigma g(S_i)$.  Conversely, given a generalized Heegaard splitting $((V_1,W_1,F_1),\cdots , (V_n,W_n,F_n))$ for $M$, one can always use the process of {\em amalgamation} to change it into a standard Heegaard splitting $(V,W,F)$ of $M$ satisfying the same equation above.  It is not necessary to recount the details of these processes for the purposes of this paper, but we will need the following fact.\\

\begin{proposition}
If $(V,W,F)$ is a weakly reducible Heegaard splitting of $M$, then $(V,W,F)$ can be untelescoped to a generalized Heegaard splitting $((V_1,W_1,F_1),\cdots , (V_n,W_n,F_n))$ such that $(V_i,W_i,F_i)$ is a strongly irreducible splitting of $M_i$ for each $1\leq i \leq n$, and the thin surfaces $S_i$ are incompressible in $M$ for each $1\leq i <n$.\\
\end{proposition}

We will call a generalized Heegaard splitting of the kind given by Proposition 2.11 {\em fully untelescoped}.\\

\begin{proposition}
If $((V_1,W_1,F_1),\cdots , (V_n,W_n,F_n))$ is a fully untelescoped generalized Heegaard splitting for $M$, and $P$ is an incompressible surface in $M$, then $P$ can be isotoped to intersect each of the thin and thick surfaces only in mutually essential simple closed curves.
\end{proposition}

\section{Primitive Heirarchies}

The propositions proved in this section amount roughly to the following cutting and pasting lemma for compression bodies:  Given incompressible surfaces with boundary, $P$ and $P'$, on the positive boundaries of compression bodies $V$ and $V'$ respectively, and an orientation preserving homeomorphism $h:P\rightarrow P'$, the resulting manifold $W=V\cup_h V'$ is a compression body if and only if $P$, considered as a surface embedded in $W$, admits a sequence of $\partial$-compressions reducing it to a collection of disks.  However, our construction of doppelg\"{a}nger surfaces in Section 4  depends on the extra information contained in a ``primitive heirarchy'' as defined below, and so this cutting and pasting lemma must likewise be made in the language of primitive heirarchies in order to ensure that our doppelg\"{a}ngers do the job we set out for them.\\

\begin{definition}
Let $P$ be a surface with non-empty boundary.

\begin{itemize}

\item A {\em multi-arc} is a finite, disjoint union $\talpha=\alpha_1\cup \cdots \cup \alpha_n$ of essential arcs properly embedded in $P$.

\item A {\em labeled multi-arc} $(\talpha,l)$ is a multi-arc with an assigned {\em labeling function} $l:\hat{\talpha}\rightarrow  \N$.

\item A {\em partitioned multi-arc} is an equivalence class of labeled multi-arcs $[\talpha, l]$ under the relation defined as follows: $(\talpha, l)\sim(\talpha', l')$ if and only if $\talpha=\talpha'$ and $l=f\circ l'$ for some bijection $f: \N\rightarrow \N$.\\

\end{itemize}
\end{definition}

\begin{definition}

Let $[\talpha, l]$ be a partitioned multi-arc on the surface $P$, and suppose that $\beta$ is an arc properly embedded in $P$ that is disjoint from $\talpha$.  Let $\iota : E(\beta, P)\rightarrow P$ be the inclusion map, let $\tilde{\gamma}$ be the multi-arc on $E(\beta, P)$ obtained from the non-boundary parallel components of $\iota^{-1}(\talpha)$, and let $\iota_{\#}:\hat{\tilde{\gamma}}\rightarrow \hat{\talpha}$ be given by $\alpha\mapsto \iota(\alpha)$.  Then the {\em projection} of $[\talpha, l]$ onto $E(\beta, P)$ along $\beta$ is defined by the equation $p_\beta ([\talpha, l])=[\tilde{\gamma}, l\circ \iota_\# ]$.  For an arc $\alpha\in \hat{\talpha}$, we define $p_\beta(\alpha)=\iota^{-1}(\alpha)$ when $\iota^{-1}(\alpha)$ is essential, and $p_\beta(\alpha)=0$ otherwise.\\

\end{definition}

\begin{definition}
If $[\talpha, l]$ be a partitioned multi-arc on a surface $P$, then for any $\alpha\in \hat{\talpha}$, let $s(\alpha)=\cup (l^{-1}\circ l (\alpha))$, and call the components of $s(\alpha)\setminus \alpha $ the {\em sister arcs} of $\alpha$.  If $\alpha$ has no sister arcs, it is said to be {\em primitive}.  If $\talpha$ itself consists of only one arc, we say that $[\talpha, l]$ is primitive.\\
\end{definition}

\begin{definition}
Let $[\talpha, l]$ be a partitioned multi-arc on a surface $P$.  A {\em primitive heirarchy} on $P$ supported by $[\talpha, l]$ is a sequence $(\alpha_1, \cdots , \alpha_n)$ of elements from $\hat{\talpha }$ satisfying the following conditions:

\begin{itemize}

\item  $\alpha_1$ is primitive in $[\talpha, l]$.

\item  For $k >1$, $p_{\alpha_{k-1}}\circ p_{\alpha_{k-2}}\circ \cdots \circ p_{\alpha_1}(\alpha_k)\neq 0$, but $p_{\alpha_{k-1}}\circ p_{\alpha_{k-2}}\circ \cdots \circ p_{\alpha_1}(\alpha)=0$ for every sister arc of $\alpha_k$ in $[\talpha, l]$.

\item  Every component of $E(\alpha_1\cup \cdots \cup \alpha_n, P)$ is a disk.\\

\end{itemize}

\end{definition}

\begin{figure}
\centering
\includegraphics[width=0.8\textwidth]{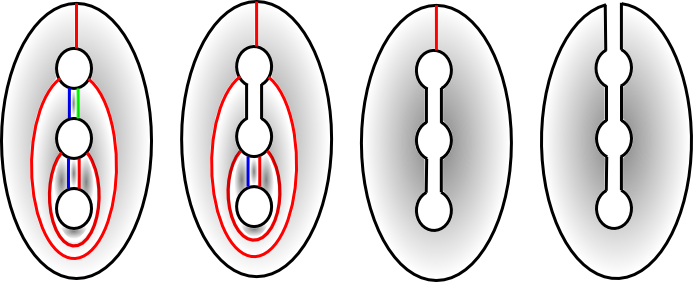}
\caption{A partitioned multi arc which supports a primitive hierarchy, as displayed in a sequence of projections.  Arcs of the same color are understood to share a common label.}

\end{figure}

It is an immediate consequence of this definition that if $(\alpha_1, \cdots , \alpha_n)$ is a primitive heirarchy on $P$ supported by $[\talpha, l]$, and $\iota: E(\alpha_1,P)\rightarrow P$ is the inclusion map, then $(\iota^{-1}(\alpha_2), \cdots , \iota^{-1}(\alpha_n))$ is a primitive heirarchy on $E(\alpha_1, P)$ supported by $p_{\alpha_1}([\talpha, l])$.

\begin{definition}

Given two partioned multi-arcs $[\talpha, l]$, $[\talpha',l']$ on a surface $P$ such that $\talpha\cap \talpha'=\emptyset$, with $im(l)\cap im(l')=\emptyset$ (which can always be ensured by making an appropriate choice of $l$ and $l'$), the {\em amalgamation} of $[\talpha,l]$ and $[\talpha',l']$ is the partitioned multi-arc $[\talpha\cup \talpha',l\cup l']$.  The amalgamation of $k$ partitioned multi-arcs on a surface is defined similarly.\\

\end{definition}

\begin{definition}
Let $P$ be a subsurface of the boundary of a $3$-manifold $M$ and let $F=F_1\cup \cdots \cup F_n$ be a disjoint union of connected surfaces $F_i$ embedded in $M$ which meet $P$ only in essential arcs on $P$.  Then if $\talpha_i=F_i\cap P$ and $l_i:\hat{\talpha}_i\rightarrow 
N$ is constant, then the amalgamation $[\bigcup_{i=1}^n\talpha_i,\bigcup_{i=1}^nl_i]$ is called the {\em partitioned multi-arc induced by $F$}.\\

\end{definition}

\begin{definition}
Given a partioned multi-arc $[\talpha, l]$ on $P$ and a homeomorphism $h:P\rightarrow P'$, the {\em pushforward of $[\talpha,l]$ onto $P'$ along $h$}, denoted $h_*[\talpha,l]$ is the partioned multi-arc $[h(\talpha),l\circ h_\#^{-1}]$, where $h_\#$ is the induced map on sets.\\
\end{definition}

The machinery defined above is recursive in nature and is thus particularly amenable to inductive arguments.  In the closing propositions of this section, we shall induct on the following surface complexity.  It is appropriate for our purposes because it decreases after compressions and $\partial$-compressions, while ignoring disks and spheres.\\

\begin{definition}
Given a surface $P$, let $\tilde{P}\subset P$ be the union of the components of $P$ which are not disks or spheres.  Define the complexity $c(P)$ of $P$ to be the pair $(-\chi (\tilde{P}), |\tilde{P}|)$, and order complexities lexicographically.\\
\end{definition}

The following proposition is well-known, although in the literature it is typically asserted in less formal language by saying that $P$ can be $\partial$-compressed down to a collection of disks.\\

\begin{proposition}
Let $P$ be an incompressible surface properly embedded in a compression body $V$ such that $\partial P\subset \partial_+V$, and each of whose connected components $P_1,\cdots P_k$ has non-empty boundary.  Let $\pi_j:Fr(P_j)\rightarrow P_j$ be the bicollar projection for $1\leq j\leq k$.  Then $E(P,V)=V_1\cup \cdots \cup V_n$, where the union is disjoint and each $V_i$ is a compression body, and for $1\leq i\leq n$ there are disjoint unions of compressing disks $\D_i$ for $V_i$ with the following properties.\\

\begin{itemize}
\item $\D_i\cap Fr(P)$ consists entirely of essential arcs.

\item $\pi_j(\D_i)\cap \pi_j(\D_m)=\emptyset$ for all $j$ and $i\neq m$.

\item If $[\talpha_i,l_i]$ is the partioned multi-arc induced by $\D_i$, and $\pi=\pi_1\cup \cdots \cup \pi_k$ then the amalgamation $[\talpha,l]$ of all the partitioned multi-arcs $\pi_*[\talpha_i,l_i]$ supports a primitive heirarchy on $P$.\\

\end{itemize}

\begin{proof}
Let $\D$ be a complete collection of disks for $V$ which intersects $P$ only in essential arcs.  Then for each component $V_i$ of $E(\D,V)$, the components of $P\cap V_i$ are incompressible with nonempty boundary.  In particular, if $V_i$ is a ball then $P\cap V_i$ consists of disks.  If $V_i$ is a thickened surface $F\times I$ parameterized so that $F\times \{0\}\subset \partial_-V$, then by Proposition 2.7 each component of $P\cap V_i$ is $\partial$-compressible.\\

We induct on $c(P\cap E(\D,V))$.  In the base case, every component of $P\cap E(\D,V)$ is a disk.  If we let $\D_i=\D\cap V_i$, then the first condition is satisfied by our hypothesis on $\D$, and after a proper isotopy of these disks supported on a small open neighborhood of $(\bigcup_{i=1}^n\D_i)\cap Fr(P)$ in $E(P,V)$ we can ensure that the second condition is satisfied as well, and so (after this isotopy) $\D'=(\bigcup_{i=1}^n\D_i)$ imposes a partitioned multi-arc on $P$.\\  

For the third condition, let $Y=P\cap \D$ and let $\Cf_Y$ be the connectivity graph on $Y$ in $\D$, which has trees as components, one for each disk component of $\D$ that intersects $P$.  The vertices of $\Cf_Y$ correspond to components of $\D'$, and given a component $D$ of $\D'$, we let $d(D)$ denote the simplicial distance of its corresponding vertex in $\Cf_Y$ from the set of leaves (so $d(D)=0$ means that $D$ is an outermost disk, and more generally $d(D)=n-1$ means that $D$ can be described as $n$-th outermost in $\D$ with respect to $Y$).  Let $\sigma:\hat{\D}'\rightarrow \N$ be an injective function which respects the ordering imposed by $d$, so that $d(D)<d(D')$ implies $\sigma(D)<\sigma(D')$.  For any $D\in \hat{\D}'$, $\pi(D\cap Fr(P))$ has at most one arc component which is not isotopic to an arc of $\pi((\bigcup_{\sigma(D')<\sigma(D)}D')\cap Fr(P))$ in $P$, call this arc $\alpha_{D}$ if it exists (see Figure 2).  Then the sequence $(\alpha_1,\cdots ,\alpha_r)$ on the set of components of $\alpha_{\D'}=\bigcup_{D\in \D'}\alpha_{D}$ satisfying $\sigma(\alpha_s)<\sigma(\alpha_t)$ whenever $s<t$ satisfies the first two conditions of Definition 3.4 required to be a primitive heirarchy on $[\talpha,l]$. Moreover, since every arc of $\D\cap P$ is properly isotopic to an arc in $\alpha_{\D'}$ in $P$, $E(\alpha_{\D'},P)$ consists entirely of disks.  Thus the base case is proved.\\

Now consider a surface $P$ such that $P\cap E(\D,V)$ does not consist entirely of disks, and assume that every surface of lower complexity satisfies the conclusions of the proposition.  Then there is an outermost $\partial$-compressing disk of $P$, and the surface $P'$ obtained by $\partial$-compressing $P$ along $D$ has lower complexity than $P$.  Thus $P'$ admits a collection of disks $\D'$ satisfying the conclusions of the proposition.  Moreover, elementary general position arguments ensure that $\D'$ can be chosen so that it intersects $N(D)$ only in disks whose boundary lies in $\partial N(D)\setminus Fr(D)$.  It is then easy to verify that the union of disks $(\D'\cap E(P,V))\cup D$ also satisfies the conclusions of the theorem.\\

\end{proof}

\end{proposition}

\begin{figure}
\centering
\includegraphics[width=0.8\textwidth]{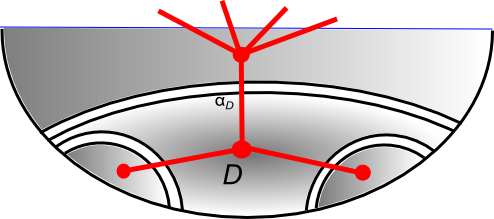}
\caption{A portion of a disk in $\D$ with the corresponding part of the connectivity graph $\Cf_Y$ superimposed in red.  The middle region corresponds to the component $D$ of $\D'$, with its edge $\alpha_D$ labeled.  In this picture $d(D)=1$.}

\end{figure}

The following becomes the converse to Proposition 3.9 after making the additional assumption that $P$ and $P'$ are incompressible.\\

\begin{proposition}
Let $V$ and $V'$ be unions of compression bodies, let $P\subset \partial_+ V$ and $P'\subset \partial_+ V'$ be subsurfaces, and let $h:P\rightarrow P'$ be a homeomorphism.  Suppose also that $V$ and $V'$ admit disjoint collections of compressing disks $\D=D_1\cup \cdots \cup D_k$ and $\D'=D_{k+1}\cup \cdots \cup D_n$, respectively, which satisfy the following properties:

\begin{itemize}
\item $h(\D\cap P)\cap (\D'\cap P')=\emptyset$, and $h(\D\cap P)\cup (\D'\cap P')$ consists only of essential arcs in $P'$.

\item The partitioned multi-arc $[\talpha, l]$ on $h(P)=P'$ induced by $\D\cup \D'$ in $V\cup_h V'$  supports a primitive heirarchy.

\end{itemize}

Then $V\cup_h V'$ is also a union of compression bodies.

\begin{proof}
W induct on the complexity $c(P)$.  In the base case every component of $P$ is a disk and the proposition is trivial.\\

So suppose $P$ is a surface with at least one non-disk component and assume the proposition holds for all surfaces of complexity lower than $c(P)$.  Let $\alpha$ be the first element of the primitive heirarchy on $P'=P\setminus \{disks\}$ supported by $[\talpha,l]$.  Then without loss of generality we may assume that there is a disk component $D$ of $\D'$ such that $D\cap P'=\alpha$.  Let $W=V\cup_hN(D)$, $W'=E(D, V')$, $F=E(h^{-1}(\alpha),P)\cup_hFr(D,V')\subset W$, and $F'=E(\alpha,P')\cup Fr(D,V')\subset W'$.  Then $W$ and $W'$ are both unions of compression bodies, and $F$ is a surface of lower complexity than $P$.\\

All boundary parallel arcs of $\D \cap F$ in $F$ can be removed via a proper isotopy of $\D$ in $W$ that leaves all the remaining arcs of $\D\cap F$ unchanged, and likewise for the boundary parallel arcs of $(\D'\setminus D)\cap F'$ in $W'$.  The homeomorphism $h|_{P\cap F}:P\cap F \rightarrow P'\cap F'$ can be extended to a homeomorphism $g: F\rightarrow F'$.  Then $g$ and the isotoped version of $\D\cup\D'\setminus D$ define a partitioned multi-arc $[\talpha', l']$ on $F'$, and after an appropriate identification of $F'$ with $E(\alpha, P')$ (one given by a homeomorphism which restricts to the identity outside of a small open neighborhood of $Fr(\alpha,P')$ in $P'\cap F'$), we see that in fact $[\talpha',l']=p_\alpha([\talpha,l])$.\\

By the remark immediately following Definition 3.4, it follows that $[\talpha',l']$ supports a primitive heirarchy on $F'$.  Since $F$ is of lower complexity than $P$, the induction hypothesis tells us that $W\cup_gW'$ is a union handlebodies.  Since $W\cup_gW'\cong V\cup_hV'$ the proof is finished.\\

\end{proof}

\end{proposition}

\section{Doppelg\"{a}ngers}

In this section, $M_1$ and $M_2$ are compact manifolds and, for $i=1,2$, $T_i\subset \partial M_i$ are incompressible tori.  Also $h: T_1\rightarrow T_2$ is a homeomorphism, $M=M_1\cup_h M_2$ with quotient map $\pi: M_1\cup M_2\rightarrow M$, and $T=\pi(T_1)=\pi(T_2)$ is the resulting incompressible torus in $M$.  In this section we prove that $g(M)\geq g(M_1)-1$ under the assumption that $M$ admits a minimal genus splitting $(F,V_1,V_2)$ satisfying the following conditions:\\

\noindent {\bf Property A:} $F$ intersects $T$ only in essential simple closed curves.  \\

\noindent {\bf Property B:} The connectivity graph $\Cf_{F_2}$ of $F_2=F\cap M_2$ in $M_2$ is simply connected.  \\

It follows from Properties A and B that $P_i=T\cap V_i$ is a union of incompressible annuli, and each component of $E(P_i,V_i)$ is a compression body.  Furthermore, by Proposition 3.9, there are collections of compressing disks $\D_i\subset V_i\cap M_1$ and $\E_i\subset V_i\cap M_2$ which induce partitioned multi-arcs $[\talpha^i,l_i]$ and $[\tbeta^i,p_i]$ on $P_i$, the almagamation of which supports primitive heirarchies $(\gamma_1^i,\cdots , \gamma_{n}^i)$ on $P_i$, for $i=1,2$.  Let $\lambda\subset T$ be an oriented simple closed curve with basepoint $x$ which intersects each component of $P_1$ and $P_2$ in a single essential arc.  Then $\lambda$ induces an indexing function $\lambda_\# :\hat{P}_1\cup \hat{P}_2 \rightarrow \N$ given by letting $\lambda_\#(P')$ be the order in which $\lambda$ intersects the component $P'$, starting from $x$.  Similarly, $\lambda$ induces a circuit $c_\lambda$ on the connectivity graph $\Cf_{F_2}$ with basepoint on the vertex corresponding to the component of $E(F_2,M_2)$ on which $x$ lies.\\

\begin{figure}
\centering
\includegraphics[width=0.8\textwidth]{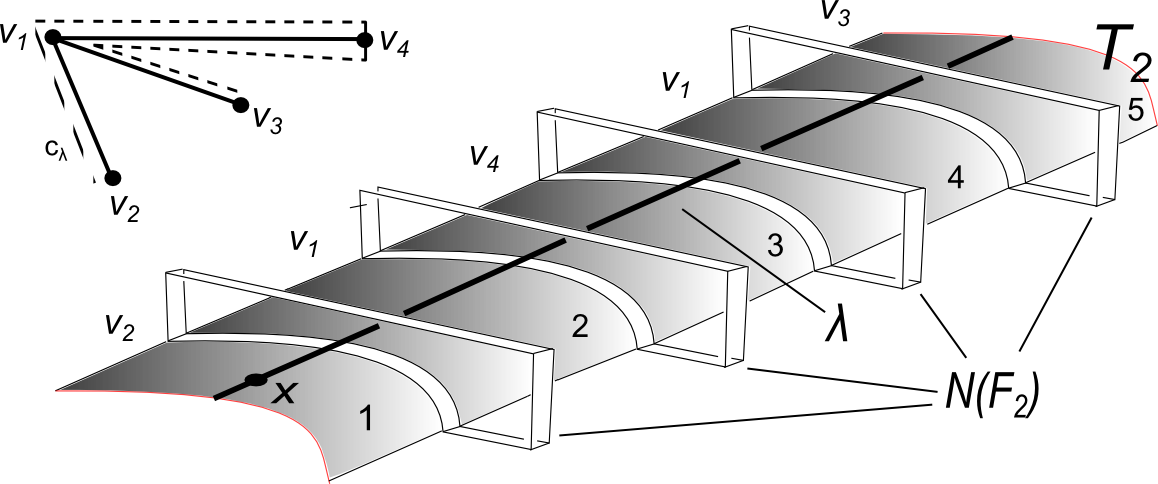}
\caption{A small region of $N(T_2)$ in $M_2$.  $T_2$ itself is shaded, $N(F_2)$ is translucent.  The annular components $P_1\cup P_2$ of $E(\partial N(F_2),T_2)$ are numbered as by $\lambda_{\#}$, and components of $E(F_2, M_2)$ are labeled with the same label as their corresponding vertex in the connectivity graph for $F_2$ in $M_2$, a portion of which appears in the upper left hand corner as determined by the diagram. The dashed line indicates the corresponding portion of $c_\lambda$ as determined by what can be seen of $\lambda$.}

\end{figure}

To prove the main proposition of this section, we construct an $(M_2,F_2)$-{\em doppelg\"{a}nger pair} $(W,Q)$ where $W$ is a solid torus and $Q$ is a properly embedded surface in $W$ with the following properties:\\

\begin{enumerate}
\item $Q\cap \partial W=\partial Q$ consists of essential closed curves on $\partial W$, and every component of $E(Q,W)$ is handlebody.\\

\item For $i=1,2$, there exists a collection of compressing disks $\D'_i$ embedded in $E(Q,W)$ which induces (after pushing forward along $f$) a partitioned multi-arc $[\tomega^i, q]$ on $P_i$ whose amalgamation with $[\talpha^i, l]$ forms a partitioned multi-arc supporting a primitive heirarchy on $P_i$.\\

\item There is an isomorphism $g:\Cf_Q\rightarrow \Cf_{F_2}$ between the connectivity graphs of $Q$ in $W$ and $F_2$ in $M_2$.  Moreover, the Euler characteristic of the surface corresponding to an edge $e$ in $\Cf_Q$ is always higher than that of the surface corresponding to $g(e)$ in $\Cf_{F_2}$.\\

\item There exists a homeomorphism $f:(T_1, \partial P_1) \rightarrow (\partial W, Q)$ which preserves the circuit $c_\lambda$ with respect to $g$.  In other words, the circuit $g(c_\lambda)$ is the same as the circuit imposed by the oriented simple closed curve $f^{-1}(\lambda)$ with basepoint $f^{-1}(x)$ on $\Cf_Q$.\\

\item The core of $W$ can be embedded in $Q$.\\

\end{enumerate}

From (1), (2) and Proposition 3.10 we deduce that every component of $E(F_1,M_1)\cup_f E(Q,W)$ is a compression body.  From (3) and (4) it follows that $F_1\cup_f Q$ is connected and of genus less than or equal to $F$.  These facts together tell us that $g(M)\geq g(M_1\cup_f W)$. After removing the core of $W$ from $M_1\cup W$ we obtain $M_1$ again, and under this identification we deduce from (5) that, after at most one stabilization, $F_1\cup_f Q$ is also a Heegaard surface for $M_1$.  It then follows that $g(M)\geq g(M_1)-1$ as claimed.  From now on we will refer to properties (1)-(5) as the {\em $Q$ properties}.\\

The remainder of the section is devoted to the construction of $(W,Q)$ and the proof that it satisfies the $Q$ properties.  We fix a parameterization of $W$ as $S^1\times \Df$, where $S^1$ is the unit circle parameterized by $0\leq\varphi<2\pi$, and $\Df$ is the closed unit disk in $\R^2$ centered at the origin, parameterized by polar coordinates $(r,\theta )$. Our gluing homeomorphism $f:\partial T_1\rightarrow \partial W$ is the unique one (up to isotopy) which sends $\lambda$ to a meridian $\{pt\}\times \partial\Df$ of $W$, and which sends the curves of $F_1\cap T_1$ to longitudes $S^1\times \{pt\}$.\\

Let $X\subset W$ be an embedded graph with one {\em central} vertex $c$ at $\{\varphi_0\}\times (0,0)\in S^1\times \Df$, a finite number $n>1$ of {\em outer} vertices $x_1,\cdots , x_n$ on $\{\varphi_0\}\times \partial D$, one {\em radial} edge connecting each outer vertex to the central vertex, and for each outer vertex $x_i$ with coordinates $(\varphi_0,1,\theta_i)$, one {\em longitudinal} edge $S^1\times (1,\theta_i)$ connecting $x_i$ to itself.  Then $X$ is said to be a {\em connected spoked graph} in $W$, and (any choice of) $Fr(X,W)$ is called a {\em connected spoked surface} for $V$ with {\em spine} $X$.  A finite, disjoint union of connected spoked graphs $\tilde{X}$ will simply be called a {\em spoked graph}, and $Fr(\tilde{X})$ then is said to be a {\em spoked surface}.  Moreover, we call $E(\tilde{X},W)$ a {\em spoked chamber}.\\

If $S$ is a connected spoked surface in $W$, then the two components of $E(S,W)$ are homeomorphic.  In fact, if $S=\partial N(X)$ for a connected spoked graph $X$, then after translating $X$ in the $\varphi$ direction and rotating it slightly in the $\theta$ direction so that the resulting graph $X'$ is disjoint from $N(X)$ and has one longitudinal edge in each annulus component of $\partial W\setminus N(X)$, $X'$ is also a spine of $S$.  We will say that $X'$ is {\em dual} to $X$.\\

It follows that if $\Sf$ is a (not necessarily connected) spoked surface with spine $\tilde{X}$, and $S$ is a component with spine $X\subset \tilde{X}$, then we can perform an ambient isotopy of $W$, which takes all of $E(X,W)$ onto a small closed neighborhood $N(X',W)$ of some dual spine $X'$ of $X$.  Then given any other spoked graph $\tilde{X}'$ disjoint from $N(X',W)$, $\Sf'=S\cup \partial N(\tilde{X'})$ will become another spoked surface embedded in $W$ having $S$ as a component.  In this way we see that $S$ can be made to cobound spoked chambers of any kind we choose on either side (subject only to the condition that they are cobounded by $S$, of course).  This process can then be repeated with each of the components of $\partial N(\tilde{X'})$, and so on again.  We refer to this as the {\em free nesting property} of spoked surfaces.\\

It is an immediate consequence of the free nesting property that, for any tree $\Tf$, we can construct a surface $Q' \subset W$ such that each component of $Q'$ is a spoked surface, each component of $E(Q',W)$ is a spoked chamber (and thus a handlebody), and whose connectivity graph $\Cf_Q'$ is isomorphic to $\Tf$.  It is clear that the $Q$ properties (1) and (5) hold for any such surface, but the remaining $Q$ properties ask for a bit more precision.\\

Let $\mu=\{0\}\times \partial \Df$, we may assume $f(\mu)=\lambda$.  The points of $f^{-1}(\partial F_1\cap \lambda)$ subdivide $\mu$ into arcs $a_1,\cdots a_{2n}$, where we assume the arc $f(a_j)$ lies on the annulus $A\in \hat{P}_1\cup \hat{P}_2$ satisfying $\lambda_\#(A)=j$.  In other words we have ordered the arcs $a_j$ to correspond to the ordering $\lambda$ imposed on the annuli $\hat{P}_1\cup \hat{P}_2$.  Finally, to any $a_j$ we associate the vertex $v(a_j)$ of $\Cf_{F_2}$ corresponding to the component of $E(F_2,M_2)$ which $\pi\circ f(a_j)$ intersects.\\

Now, given an arbitrary vertex $v$ of $\Cf_{F_2}$, we construct a spoked chamber in $W$ as follows.  First let $\Af=\{a_{j_1},\cdots , a_{j_l}\}$ be the set of all subarcs satisfying $v(a_j)=v$.  Since $\Cf_{F_2}$ is a tree, as we travel around $\mu$ from $a_{j_s}$ to $a_{j_{s+1}}$, we know that $v(a_{{j_s}+1})=v(a_{{j_{s+1}}-1})$.  That is, $\mu$ always comes back the same way it came.  Thus we can unambiguously associate the vertex $v(a_{{j_s}+1})$ of $\Cf_{F_2}$ to the complementary subarc of $\cup \Af$ in $\mu$ which lies between $a_{j_s}$ and $a_{j_{s+1}}$.  The same can be done for all of the complementary subarcs $r_1, \cdots r_l$ of $\cup \Af$ in $\mu$, and we let $v(r_j)$ denote the vertex of $\Cf_{F_2}$ associated to $r_j$ in this fashion.  We let $\Wf$ be the set of vertices in $\Cf_{F_2}$ adjacent to $v$.\\

Let $y_j$ be a longitude of $\partial W$ lying in $S^1\times r_j$.  For each subset of $R=\{r_1,\cdots , r_l\}$ of the form $R_w=\{r\in R: v(r)=w\}$, $w\in \Wf$, we can embed a connected spoked graph $X_w$ in $W$ which has central vertex $c_w=(\varphi_w,0,0)$ and longitudinal edges $y_j$ associated with the set $R_w$.  Moreover we choose $N(X_v)$ so that $N(X_w)\cap \partial W=S^1\times \overline{\cup R_w}$, which ensures that the intersection of $\partial N(X_w)$ with $\partial W$ lies in $f^{-1}(F_1\cap T_1)$.  We now need to specify the values of $\varphi_w$ more closely in order to ensure that $Q$-property $2$ will hold.\\

Now, for every $a_j$, $1\leq j\leq 2n$ we can let $\gamma(a_j)$ be the arc $\gamma_k^i$ which lies on the same element of $\hat{P}_1\cup \hat{P}_2$ as $f(a_j)$, for $i=1,2$ (recall that $(\gamma_1^i,\cdots , \gamma_1^n)$ was our primitive hierarchy on $P_i$).  This correspondence is a bijection.  Moreover, since $v(a_{j_s})=v$ for all $a_{j_s}\in \Af$, we may assume without loss of generality that $\gamma(a_{j_s})$ is always one of the arcs $\gamma_k^1$.  Call an arc $a_{j_s}$ a {\em connector} if the subarcs $r_s$, $r_{s+1}$ adjacent to it have the property $v(r_s)\neq v(r_{s+1})$, and in this case say that $a_{j_s}$ {\em connects $v(r_s)$ to $v(r_{s+1})$}.  Let $\Bf=\{b_1,\cdots , b_t\}\subset \Af$ be the set of all connectors, indexed so that $i<j$ implies that $\gamma(b_i)$ appears earlier in the primitive hierarchy $(\gamma_1^i,\cdots , \gamma_1^n)$ than $\gamma(b_j)$.\\

The connector $b_1$ will connect two distinct vertices of $\Cf_{F_2}$, call them $w_0$ and $w_1$, and let $\varphi_{w_0}=0$, $\varphi_{w_1}=\pi$.  Suppose now we have defined $\varphi_w$ for every vertex $w\in \Wf_k$, where $\Wf_k$ is the set of all vertices in $\Cf_{F_2}$ which are connected to some other vertex in $\Cf_{F_2}$ by an element in $\{b_1,\cdots, b_k\}$ (so for example $\Wf_1=\{w_0,w_1\}$).  Then there are three possibilities we must consider when defining $\varphi_w$ for the elements $w\in \Wf_{k+1}$.\\

The first possibility is that the two vertices, call them $w$ and $w'$, connected by $b_{k+1}$ both already lie in $\Wf_k$, so that $\Wf_{k+1}=\Wf_{k}$ and there is nothing to do.\\

The second possibility is that neither of $w$ or $w'$ lies in $\Wf_k$, in that case define $\varphi_w<2\pi$ and $\varphi_{w'}<2\pi$ so that both are distinct and larger than any of the values $\varphi_{w''}$ for $w''\in \Wf_k$.\\  

The third possibility is that $w$ lies in $\Wf_k$ while $w'$ does not, and here it is necessary to look a little bit closer at $\Wf_k$.  Let $\sim$ be the smallest equivalence relation on $\Wf_k\times \Wf_k$ containing the points $(w'',w''')$ where $b$ connects $w''$ to $w'''$, and call the subsets in the partition of $\Wf_k$ coming from $\sim$ {\em cluster sets}.  If $\tilde{\Wf}$ is the cluster set containing $w$, then define $\varphi_{w'}<2\pi$ to be larger than any angle $\varphi_{w''}$ for $w''\in \tilde{\Wf}$, but smaller than any other angle $\varphi_{w'''}$ having the same property, for $w'''$ in $\Wf_k$.\\

With the set of values for $\varphi_w$ are inductively defined as above for $w\in \Wf$, we will say that the spoked graph $\tilde{X}_v=\bigcup_{w\in\Wf} X_w$ is {\em well clustered} with respect to the primitive hierarchies $(\gamma_1^i,\cdots , \gamma_n^i)$, and likewise for the chamber $W_v=E(\tilde{X}_v,W)$ associated to $v$.  Using the free nesting property on each of the surfaces $\partial N(X_w)$, we can repeat the entire process for each vertex in $\Wf$, and continue in this way along all of $\Cf_{F_2}$ until we have built $W_z$ and its associated surfaces for every vertex $z$ in $\Cf_{F_2}$.  The union of the resulting surfaces will be $Q$, and it is clear that the $Q$-properties (3) and (4) hold for $Q$.  The only thing left to verify is $Q$-property (2).  For simplicity we will restrict our attention to what is happening in $V_1$, the proof for the $V_2$ side is identical.\\

I will induct on $n$, the number of components of $P_1$.  In the base case $n=1$, $Q$ is just an annulus with longitudinal boundary components which cuts $Q$ into two solid tori components $W_1,W_2$, and we can let the disk set $\D_1'$ be the meridian disk of $W_1$, isotoped so as to ensure it intersects $f^{-1}(P_1)$ in a single arc disjoint from $\talpha_1$.\\

Now assume Property (2) holds in our construction for all values less than $n$, $n>1$.  Let $\gamma_1^1$ be the first element in our primitive hierarchy, and suppose that $\gamma_1^1$ was one of the arcs $\beta$ in $\tbeta^1$.  Then $\beta$, being primitive, cobounds a disk $D$ in some component of $E(F_2,M_2)$ which is otherwise disjoint from $P_1$, and thus $\partial \beta$ lies on a single component of $F_2$.  This then implies that the boundary of $f^{-1}\circ h^{-1} (\beta)$ lies on a single component $Q'$ of $Q$.\\

Let $f^{-1}(P')$ be the component of $f^{-1}(P_1)$ on which $f^{-1}\circ h^{-1} (\beta)$ lies, and let $W'$ be the component of $E(Q,W)$ whose boundary meets $f^{-1}(P')$.  Since $Q'$ is a spoked surface and $\beta$ is an essential arc on an annulus lying between a pair of adjacent radial edges of the spine of $Q'$, we have an obvious compressing disk $D$ with $\partial D=\beta\cup \beta'$, where $\beta'=D\cap Q$.  $E(D,W')$ can then itself be regarded as a spoked chamber whose surfaces are well clustered with respect to the primitive hierarchy $(\gamma_2^1,\cdots , \gamma_n^1)$ on the surface $P_1\setminus P'$.  By the inductive hypothesis we then have a disk set $\D_1''$ satisfying Property (2) with respect to this hierarchy, which can be chosen disjoint from $D$ by standard transversality arguments (we merely need $\partial \D_1''\cap Q''$ not to intersect the disks $Fr(D)$, where $Q''$ is the surface obtained from $Q'$ via a compression along $D$, i.e. $Q''=E(\partial D, Q')\cup Fr(D)$).  It then follows that $\D_1'=D\cup \D_1''$ is a disk set satisfying Property (2) with respect to the primitive hierarchy $(\gamma_1^1,\cdots , \gamma_n^1)$ on $P_1$, as required.\\

To complete the induction, we must consider the second possibility where $\gamma_1^1$ is one of the arcs $\alpha$ in $\talpha^1$.  Let $P'$ be the component of $P_1$ on which $\alpha$ lies, and let $W'$ be the component of $E(Q,W)$ meeting $f^{-1}(\alpha)$.  If the boundary of $f^{-1}(\alpha)$ does meet the same component of $Q$, then we argue as before.  However it is possible that it meets distinct components $Q'$ and $Q''$ of $Q$, and here is where we use the fact that $Q$ is well clustered  For if $X_{Q'}$ and $X_{Q''}$ are the spines associated with $Q'$ and $Q''$ in our construction of the chamber $W'$, then the fact that $Q$ is well-clustered tells us that there are no components of $Q$ lying between $Q'$ and $Q''$ in the positive $\varphi$ direction.  Moreover, the surface $Q'''=\partial E(P',W')$, is also a spoked surface, and in fact it is well-clustered with respect to the primive hierarchy $(\gamma_2^1,\cdots , \gamma_n^1)$ of the surface $P_1\setminus P'$, and so by the inductive hypothesis we obtain a collection of disks $\D_1''$ in $E(f^{-1}(\alpha ),W')$ satisfying Property (2).  A component $D$ of $\D_1''$ will be properly embedded in $W'$ unless $\partial D$ meets $Fr(P',W')$, although we can assume this meeting consists only of essential arcs on $Fr(P',W')$.  In this case, we may extend each such disk into $N(P',W')$ to obtain a new disk properly embedded in $W'$ and disjoint from $f^{-1}(\alpha)$.  Let $\D_1'$ be the collection of disks resulting from this extension, completing the proof.\\

We have now completed our construction of the Doppelg\"{a}nger $(W,Q)$ of $(M_2,F_2)$, and have thus proved the main theorem of this section.  Moreover, since every component of $Q$ is planar and has the same number of boundary components as the corresponding component of $F_2$, we can in fact deduce $g(M)\geq g(M_1)-1+g(F_2)$, where $g(F_2)$ denotes the sum of the genera of the components of $F_2$.\\

\section{Applications}

\begin{definition}
Let $M$ be a compact $3$-manifold and $T$ a torus boundary component of $M$.  

Define the {\em $\Ef$-slopes} of $(M,T)$, denoted $\Ef(M,T)$, to be the set of {\em essential} slopes $[\alpha]$ on $T$ for which there exists a connected, orientable, incompressible, non-boundary parallel surface $F$ properly embedded in $M$ such that $\partial F\subset T$ and such that $\alpha$ can be homotoped onto each component of $\partial F$ in $T$.

Define the {\em $\Hf$-slopes} of $(M,T)$, denoted $\Hf(M,T)$, to be the set of slopes $[\alpha]$ in $T$ for which there exists a non-separating, orientable surface $F$ properly embedded in $M$ such that $\partial F\subset T$ and such that $\alpha$ can be homotoped onto each component of $\partial F$ in $T$.

\end{definition}

\begin{proposition}
Let $M$ be a compact $3$-manifold and $T$ a torus boundary component of $M$.  Then $\Hf(M,T)\subset \Ef(M,T)$.

\begin{proof}
If $F$ is a non-separating, orientable surface with $\partial F\subset T$ and boundary components of slope $[\alpha]$, then after maximally compressing $F$ and removing any separating components, we will obtain an incompressible surface $F'$ with the same properties.  Since $F'$ is non-separating, it cannot be boundary parallel.\\
\end{proof}

\end{proposition}

\begin{definition}
Let $M_1$ and $M_2$ be compact $3$-manifolds and let $T_1$ and $T_2$ be torus boundary components of $M_1$ and $M_2$, respectively.

A homeomorphism $h:T_1\rightarrow T_2$ is said to be $\Ef$-incongruous if $h_*(\Ef(M_1,T_1))\cap \Ef(M_2,T_2)=\emptyset$, where $h_*$ is the map induced on homology.\\

\end{definition}

If $M_1$ and $M_2$ are compact irreducible $3$-manifolds with $\partial M_i$ consisting entirely of tori for $i=1,2$, then a famous result of Hatcher's (reference) implies that $\Ef$-incongruous maps are, in a strong sense, generic.  In fact, when $\partial M_i$ consists only of the torus $T_i$ for $i=1,2$, Hatcher's result tells us that all but a finite number of isotopy classes of homeomorphisms $h: T_1\rightarrow T_2$ will be $\Ef$-incongruous.\\

\begin{proposition}
Let $M_1$ and $M_2$ be compact $3$-manifolds and let $T_1$ and $T_2$ be incompressible torus boundary components of $M_1$ and $M_2$, respectively.  If $h:T_1\rightarrow T_2$ is $\Ef$-incongruous, and $M=M_1\cup_h M_2$ then $g(M)\geq \min \{g(M_1),g(M_2)\}-1$.\\

\begin{proof}
Let $\pi:M_1\cup M_2\rightarrow M_1\cup_h M_2$ be the quotient map, and let $T=\pi(T_1)=\pi(T_2)$.  It is a straightforward excercise to verify that a prime factorization of $M$ can be found whose decomposing spheres $\Sf$ are disjoint from $T$ because it is incompressible.  Thus, by Haken's lemma, it suffices to prove the proposition in the case that each of $M_1$ and $M_2$ is irreducible, in which case $M$ is also irreducible, again by the incompressibility of $T$.\\

By Proposition 2.11, any minimal genus Heegaard splitting $F$ of $M$ can be fully untelescoped to a generalized Heegaard splitting $((V_1,W_1,F_1),\cdots , (V_n,W_n,F_n))$ whose thick and thin surfaces intersect $T$ only in mutually essential simple closed curves, and for one of the manifolds $M_1$ or $M_2$, the $\Ef$-incongruity of $h$ ensures us that the slope of intersection is not an $\Ef$ slope.  Suppose it is $M_1$.  Then if a thin surface intersects $T$ it intersects $M_1$ only in boundary parallel annuli, and thus after minimizing the intersection of the thick and thin surfaces with $T$ every thin surface will be disjoint from $T$.\\ 

This leaves only two possibilities.  The first is that $T$ is disjoint from every thick surface as well, in which case its incompressibility implies that it is isotopic to a thin surface, and we can deduce $g(M)\geq g(M_1)+g(M_2)-1$.  The second is that $T$ intersects a single thick surface, say $F_1$.  In this case we restrict our attention to the manifold $M'=(V_1\cup W_1)\subset M$.  Since the slope determined by $F_1\cap T$ is not in $\Ef(M_1,T)$, it is not in $\Hf(M_1,T)$, and therefore also not in $\Hf(M_1\cap M',T)$.  It follows that the connectivity graph of $F_1\cap (M_1\cap M')$ is simply connected, and thus we can construct a doppelg\"{a}nger pair $(W,Q)$ for $(M_1\cap M',F_1\cap (M_1\cap M'))$ as in Section 4, and deduce that $g(M)\geq g(M_2)-1$.\\
\end{proof}

\end{proposition}

If the connectivity graph of $F_1\cap M_i$ were known to be simply connected for $i=1,2$ in the proof above, we could have deduced $g(M)\geq \max \{g(M_1),g(M_2)\}-1$.  We now restrict our attention to cases where we have more control over the intersection $F_1\cap T$.\\

\begin{definition}
Given two manifolds $M_1$ and $M_2$ with boundary, let $A_1\subset \partial M_1$ and $A_2\subset \partial M_2$ be annuli.  Then the {\em annular sum} $(M_1,A_1)\#_\partial (M_2,A_2)$ of $(M_1,A_1)$ and $(M_2,A_2)$ is $M_1\cup_h M_2$, where $h:A_1\rightarrow A_2$ is a homeomorphism.\\
\end{definition}

Assuming as always that $M_1$ and $M_2$ are oriented and $h$ is orientation preserving, there is only one annular sum $(M_1,A_1)\#_\partial (M_2,A_2)$, which is why we omit $h$ from our notation.  If $K_1\subset M_1$ and $K_2\subset M_2$ are knots, then it is well known that $E(K_1\# K_2)$ is just the annular sum of $E(K_1,M_1)$ with $E(K_2,M_2)$ along meridional annuli.  One can similarly define a connected sum between graphs embedded in a manifold which depends on the choice of edges and orientations that reduces to an annular sum.\\

\begin{proposition}
Suppose $A$ is an annulus embedded in $M_1$ such that $A\cap M_1=\partial A\cap M_1$ is a single core curve of an annulus $A_1\subset \partial M_1$, and let $c$ be the other component of $A$.  Suppose $A_2$ is an annulus contained in a torus boundary component $T$ of $M_2$.  Then the annular sum $(M_1,A_1)\#_\partial (M_2,A_2)\cong E(c,M_1)\cup_h M_2$, where $h:\partial N(c)\rightarrow T$ is a homeomorphism mapping the curve $A\cap \partial N(c)$ to the core of $A_2$.

\begin{proof}
Choose a regular neighborhood $N(A,M_1)=V$ just large enough so that $N(c)\subset \mathring{N}(A,M_1)$ and so that $V\cap \partial M_1=A_1$.  The manifold $E(c,V)\cup_h M_2$ is homeomorphic to $M_2$ since $E(c,V)$ is just a thickened torus, and likewise $M_1$ is homeomorphic to $\overline{M_1\setminus V}$.  Since $\overline{M_1\setminus V}$ meets $(E(c,V)\cup_h M_2)$ in an annulus $A'$ which is isotopic both to $A_1$ in $M_1$, and $A_2$ in $E(c,V)\cup_h M_2$, we have $E(c,M_1)\cup_h M_2=(\overline{M_1\setminus V}, A')\#_\partial (E(c,V)\cup_h M_2,A')\cong (M_1,A_1)\#_\partial (M_2,A_2)$.\\

\end{proof}

\end{proposition}

\begin{definition}
A compact $3$-manifold $M$ is said to be a {\em near homology sphere}, or simply an $\Hf$-sphere, if every oriented, closed surface properly embedded in $M$ is separating.\\
\end{definition}

The following proposition allows us to control the slope in which generalized Heegaard surfaces intersect our gluing tori.\\

\begin{proposition}
Suppose $M_1$ is a $3$-manifold, $T_1\subset \partial M_1$ is an incompressible torus, and that there is an annulus $A$ properly embedded in $M_1$ with one boundary component $c$ on $T_1$ and the other on a different component of $\partial M_1$.  Suppose $M_2$ is another $3$-manifold with incompressible torus boundary component $T_2$, $h: T_1\rightarrow T_2$ is any homeomorphism, $\pi: M_1\cup M_2\rightarrow M_1\cup_h M_2$ is the quotient map, and $T=\pi(T_1)=\pi(T_2)$.  Then if $((V_1,W_1,F_1),\cdots , (V_n,W_n,F_n))$ is a fully untelescoped Heegaard splitting of $M$, every thick and thin surface can be isotoped to intersect $T$ only in arcs homotopic to $\pi(c)$ in $T$.

\begin{proof}
If $M'$ is the manifold obtained from $M_1$ via a Dehn filling of $M_1$ along any slope of $T_1$ whose geometric intersection number $[c]$ is $1$, then $M_1\cup_h M_2$ can instead be viewed as an annular sum between $M'$ and $M_2$ as in the previous proof.  More specifically, if we let $V=N(A\cup T_1,M_1)$, then the annulus $A'=\overline{\partial V\setminus \partial M_1}$ can be viewed as the annulus along which the annular sum is taken as in the proof of Proposition 5.6, where we have identified $M'$ with $\overline{M_1\setminus V}$.  Since $T$ is incompressible in $M$, so must be $\pi(A')$, so we may isotope the thick and thin surfaces of $((V_1,W_1,F_1),\cdots , (V_n,W_n,F_n))$ so that they intersect $\pi(A')$ only in essential simple closed curves.  Choose a collar neighborhood $W$ of $\partial V\setminus T_1$ in $V$ small enough so that $\pi(W)$ intersects the thick and thin surfaces only in vertical annuli and isotope $T$ onto $\pi(\partial W\setminus \partial V)$.  Reversing this isotopy ambiently induces an isotopy of the thick and thin surfaces which brings them to intersect $T$ in the desired fashion.\\
\end{proof}

\end{proposition}

\begin{definition}
Let $K$ be a knot in a compact $3$-manifold $M$.  $K$ is said to be $m$-small if the meridian slope of $K$ does not lie in $\Ef(E(K,M),\partial N(K))$.\\
\end{definition}

\begin{definition}
A knot $K\subset M$ is said to be {\em strongly non-trivial} if $\partial N(K)$ is incompressible in $E(K,M)$.\\
\end{definition}

\begin{proposition}
Suppose $K$ is a strongly non-trivial $m$-small knot in the $\Hf$-sphere $M_1$ with meridional annulus $A_1\subset \partial N(K)$, and $X$ is a graph in the $\Hf$-sphere $M_2$ with meridional annulus $A_2\subset \partial N(X)$.  Then $g((E(K,M_1),A_1)\#_\partial (E(X,M_2),A_2))\geq \max \{g(E(K))-1,g(E(X))\}$, and in the case that $X$ is a knot, $g((E(K,M_1),A_1)\#_\partial (E(X,M_2),A_2))\geq \max \{g(E(K)),g(E(X))\}$\\

\begin{proof}
Let $\mu\subset \mathring{E}(X,M_2)$ be the boundary of a disk $D$ that intersects $\partial N(X)$ in a core curve of $A_2$, and which intersects $X$ in a single point.  Let $[\mu']$ be the slope of $\partial N(\mu)$ determined by $D\cap \partial N(\mu)$.  We will think of $(E(K,M_1),A_1)\#_\partial (E(X,M_2),A_2)$ as $E(K,M_1)\cup_h E(X\cup \mu)$ via a gluing map $h:\partial N(K)\rightarrow \partial N(\mu)$ which sends the meridian of $K$ to $\mu'$.  Let $T$ be the image of $\partial N(\mu)$ and $\partial N(K)$ under the quotient map $\pi$ as usual.  By hypothesis $\partial N(K)$ is incompressible in $E(K,M_1)$.  Since the torus $\partial N(\mu)$ is compressible along $D$ in $E(\mu,M_2)$, any compressing disk $D'$ for $\partial N(\mu)$ must meet $\partial N(\mu)$ in the same slope $[\mu']$ as $D$.  Therefore, if $\partial N(\mu)$ were compressible along some disk $D'$ in $E(X\cup \mu, M_2)$, we obtain a sphere $S=D\cup D'\cup A'$, where $A'\subset N(\mu)$ is an annulus such that $\partial A'=\partial D\cup \partial D'$.  But $S$ then meets $X$ in a single point, which implies that $S$ is non-separating in $M_2$, contrary to our hypothesis the $M_2$ is an $\Hf$-sphere.  It follows that $\partial N(\mu)$ is incompressible in $E(X\cup \mu, M_2)$, and thus $T$ is incompressible.\\

Let $((V_1,W_1,F_1),\cdots , (V_n,W_n,F_n))$ be a fully untelescoped generalized Heegaard splitting of $(E(K,M_1),A_1)\#_\partial (E(X,M_2),A_2)$.  Then by Proposition 5.8 we may assume that its thick and thin surfaces intersect $T$ only in curves isotopic to the image of the meridian of $K$ under the quotient map, and so as in the proof of Proposition 5.4 we may isotope all the thin surfaces away from $T$.  We assume that $T$ intersects one thick surface, say $F_1$, otherwise we could deduce a stronger inequality as in the proof of Proposition 5.4.  Since $M_1$ is an $\Hf$-sphere and $F_1$ meets $\partial N(K)$ only in meridians, the connectivity graph of $F_1\cap E(K,M_1)$ is simply connected, and so we can construct a doppelg\"{a}nger pair $(W,Q)$ for $(E(K,M_1), F_1\cap E(K,M_1)$.  However, attaching the solid torus $W$ to $E(X\cup \mu, M_2)$ in the way prescribed in Section 4 amounts to a Dehn filling along $\partial N(\mu)$ whose filling slope has geometric intersection number $1$ with $[\mu']$, and this simply gives back $E(X,M_2)$.  Thus $g((E(K,M_1),A_1)\#_\partial (E(X,M_2),A_2))\geq g(E(X, M_2))$.\\

We can also construct a doppelg\"{a}nger pair $(W',Q')$ for $(E(X\cup \mu, M_2), F_1\cap E(X\cup \mu, M_2))$, because the connectivity graph of $F_1\cap E(X\cup \mu, M_2)$ must also be simply connected.  For if not, then some component, say $F$, of $F_1\cap  E(X\cup \mu, M_2)$, will be non-separating in $E(X\cup \mu, M_2)$.  If $D$ is the disk defined in the first paragraph, we may assume that $F$, considered as a surface in $E(\mu, M_2)$, intersects $D$ only in simple closed curves in its interior.  If an innermost such curve of this intersection bounds a disk $\Delta\subset D$ which does not contain the point $K\cap D$, then the result of compressing $F$ along $\Delta$ will have at least one non-separating component which we replace $F$ with.  We may continue this process until we arrive at a non-separating surface $F'$ which intersects $D$ only in curves which bound disks $\Delta \subset D$ which do contain $K\cap D$.  Let $\alpha$ be an outermost such curve of intersection, which will cobound an annulus $A'\subset D$ with $\partial D$.  The surface obtained by compressing $F'$ along $A'$ will contain at least one non-separating component, which we replace $F'$ with.  This process can be continued until we arrive at a non-separating surface $F''$ such that $D\cap F''=\emptyset$.  If two components of $\partial F''$ meet $\partial N(\mu)$ with opposite orientation induced by $F''$ and cobound an annulus $A''\subset \partial N(\mu)$ disjoint from the rest of $\partial F''$, then the surface which results from isotoping the surface $F''\cup A''$ away from $\partial N(\mu)$ near $A''$ will again be an orientable non-separating surface.  Repeating this process eventually results in a surface $F'''$.  $F'''$ cannot be closed for this would contradict the fact that $E(X\cup \mu, M_2)$ is an $\Hf$-space, thus $\partial F'''$ meets $\partial N(\mu)$ in curves all having the same orientation induced by $F'''$.  But then if we consider $F'''$ as a surface embedded (non-properly) in $E(X,M_2)$, we can create yet another orientable surface $F^4$ by attaching parallel copies of the annulus $D\cap E(X,M_2)$ along the boundary components of $F'''$.  $F^4$ is then a surface properly embedded in $E(X,M_2)$, and moreover every component of $\partial F^4\subset$ meets $\partial N(X)$ in meridional curves all having the same orientation, which implies that $F^4$ is non-separating in $E(X,M_2)$.  But then we may attach meridian disks of $N(X)$ to $F^4$ and thereby obtain a non-separating closed orientable surface $F^5\subset M_2$, which is absurd.  Thus the connectivity graph of $F_1\cap E(X\cup \mu, M_2)$ must be simply connected after all, and the proof is finished.\\

If $X$ is actually a knot, then the final statement of the proposition follows easily from Proposition 5.12 below.\\

\end{proof}

\end{proposition}

\begin{proposition}
Suppose $K_1$ and $K_2$ are non-trivial knots in $\Hf$-spheres $M_1$ and $M_2$, respectively, and that there exists a minimal genus Heegaard surface $F$ of $E(K_1\# K_2, M_1\# M_2)$ such that $F\cap A$ consists only of essential curves, where $A$ is the decomposing annulus of the connected sum.  Then $g(E(K_1\# K_2, M_1\# M_2)\geq \max \{g(E(K_1, M_1), g(E(K_2,M_2))\}$.\\

\begin{proof}
Let $E(K_1\# K_2, M_1\# M_2)$ be obtained as $E(K_1\cup \mu)\cup_h E(K_2)$ via a gluing $h: \partial N(\mu)\rightarrow \partial N(K_2)$ for some meridian $\mu$ of $K_1$ as usual, let $\pi$ be the quotient map as usual, and let $T=\pi(\partial N(\mu))=\pi(\partial N(K_2))$ as usual.  Then, as in the proof of Proposition 5.8, we may isotope $T$ so that it intersects $F$ only in essential simple closed curves with a desirable slope which ensures that the connectivity graphs of $F\cap E(K_1\cup \mu)$ and $F\cap E(K_2)$ are both simply connected, and so as in Proposition 5.11 we obtain $g(E(K_1\# K_2))\geq g(E(K_1))$.  Reverse the roles of $K_1$ and $K_2$ in the argument above to complete the proof.\\
\end{proof}

\end{proposition}

The following observation is an immediate consequence of Propositions 2.12 and 5.12.

\begin{proposition}

If $K_1$ and $K_2$ are knots in $\Hf$-spheres $M_1$ and $M_2$, and $g(E(K_1\# K_2))< \max (g(E(K_1), g(E(K_2)))$, then every minimal genus Heegaard splitting of $E(K_1\# K_2)$ is weakly reducible and intersects the decomposing annulus inessentially.\\

\end{proposition}

\section{Concluding Observations}

We have been careful to stay nearly as general as possible throughout this paper.  We have not, for example, required our manifolds be closed, so that all of the results above on connected sums of knots and graphs can be applied equally well to connected sums of links.  Also the class of $\Hf$-spheres is fairly large.  To begin with, any compact $3$-manifold $N$ embedded in a compact $3$-manifold $M$ satisfying $H_2(M,\partial M)=0$ is an $\Hf$-sphere.  We also have the following proposition.

\begin{proposition}
If $M$ is an $\Hf$-sphere with a torus boundary component $T$, then for all but at most one slope on $T$, a Dehn filling of $M$ along $T$ yields another $\Hf$-sphere.

\begin{proof}
First note that if $F\subset M$ is an orientable, non-separating surface properly embedded in $M$ with $\partial F\subset T$, then $[\partial F]\neq 0 \in H_1(T)$.  For if $[\partial F]=0$ then we can find two components $\alpha$ and $\beta$ of $\partial F$ with opposite $F$-induced orientations on $T$ and which cobound an annulus $A\subset T$ disjoint from the remainder of $\partial F$.  The surface $F'$ obtained by pushing the surface $F\cup A$ away from $T$ near $A$ is then also an orientable and non-separating surface with $[\partial F']=0$ in $H_1(T)$.  Repeating this process eventually yields a closed, orientable, non-separating surface in $M$, which is absurd.\\

The ``Half lives, half dies'' lemma tells us that for any compact orientable $3$-manifold $M'$, the image of the homology boundary map $H_2(M',\partial M')\rightarrow H_1(\partial M')$ is half the rank of $H_1(\partial M)$ (reference).  Thus there is only one possible slope $[\alpha]$ in which the boundary components  of an orientable, non-separating surface embedded in $M$ with $\partial F \subset T$ can intersect $T$.\\

So suppose a Dehn filling of $M$ along a slope in $T$ different from $[\alpha]$ yields a manifold $M'$ which is not an $\Hf$-sphere, and let $F$ be a closed, orientable, non-separating surface in $M'$.  If $V$ is the surgery torus, then by shrinking $V$ if necessary we can assume it meets $F$ only in meridian disks.  But this means that $F'=F\cap M$ is non-separating even though $\partial F'$ intersects $T$ in a slope different from $[\alpha]$, if at all.  This is absurd.\\

\end{proof}

\end{proposition}

There is good reason to believe that further progress can be made on the subject considered here, and using similar methods to those employed here.  In particular, the following proposition can very likely be proven true via a strengthening of the doppelg\"{a}nger construction of Section 4.\\

\begin{conjecture}
Given any pair of non-trivial knots $K_1$ and $K_2$ in $\Hf$-spheres $K_1$ and $K_2$, $g(E(K_1\# K_2))\geq \max \{g(E(K_1)),g(E(K_2))\}$.\\
\end{conjecture}

Conjecture 6.2 is not true if we drop the condition that $K_1$ and $K_2$ lie in $\Hf$-spheres, as is clear by setting $K_1=S^1\times \{pt\}\subset S^1\times S^2$, nevertheless it may be possible to substitute a weaker condition (notice, for example, that $\partial N(S^1\times \{pt\})$ is compressible in $S^1\times S^2$).  In any event, the very first instance of tunnel number subadditivity, found by Morimoto in (reference), shows that this conjecture is best possible, in fact it even shows that the lower bound of Proposition 5.11 is best possible.  Moreover, the subadditivity examples of (reference, reference, and reference) all realize this lower bound and allow for, or even require, that one of the knots involved is $m$-small.\\

With much less evidence or hope of a proof, we say that the following much stronger proposition might also be true.\\

\begin{conjecture}
Let $M_1$ be a compact orientable $3$-manifold with incompressible boundary component $F_1$, let $M_2$ be a compact orientable $3$-manifold with incompressible boundary component $F_2$, and let $h:F_1\rightarrow F_2$ be any homeomorphism.  Then $g(M_1\cup_h M_2)\geq \max \{g(M_1),g(M_2)\}-g(F_1)$.\\
\end{conjecture}

It is more plausible that Conjecture 6.3 can be proven in the case $g(F_1)=1$ using augmented versions of the doppelg\"{a}nger technique employed in this paper, and it is highly plausible that it can be proven in the case $g(F_1)=1$ over a suitably restricted class of compact orientable $3$-manifolds and gluing maps.  If true, Conjecture 6.3 is best possible, as has already been shown by the examples of Schultens and Wiedmann.  In fact there are even more examples.\\

\begin{proposition}
There exists an infinite family of knot pairs $K_1$ and $K_2$ in $S^3$ and gluing maps $h:\partial N(K_1)\rightarrow \partial N(K_2)$ such that $g(M)=g(E(K_1))-1$, where $M=E(K_1)\cup_h E(K_2)$.\\

\begin{proof}
It is a result of Nogueira \cite{nog} that there is an infinite family of knots $K_1\subset S^3$ with the following properties:\\

\begin{enumerate}
\item There exists a sphere $S_1\subset S^3$ intersecting $K_1$ exactly four times, and such that for each component $B$ of $E(S_1,S^3)$, $E(K_1\cap B, B)$ is a handlebody.

\item $g(E(K_1))=4$\\

\end{enumerate}

Let $K_2$ be any $2$-bridge knot in $S^3$ with $2$-bridge sphere $S_2$, and let $h:\partial N(K_1)\rightarrow \partial N(K_2)$ be any homeomorphism such that $h(S_1\cap \partial N(K_1))=S_2\cap \partial N(K_2)$.  Then, for each component $B'$ of $E(S_2, E(K_1))$, the annuli $Fr(K_2\cap B')$ are primitive because $S_2$ is a bridge sphere.  By Proposition 3.10, it then follows that the surface $(S_1\cap E(K_1))\cup_h (S_2\cap E(K_2))$ is a Heegaard surface for $M$, one whose genus is only $3$.\\

\end{proof}

\end{proposition}

It is also worth remarking that there are cases where a manifold's Heegaard genus repeatedly decreases under gluings along its incompressible boundary tori.\\

\begin{proposition}

For arbitrary $n>0$ there exists an infinite family of $n$-component links $L=L_1\cup \cdots \cup L_n\subset S^3$, $n$ knots $K_1,\cdots K_n$ in $S^3$, and homeomorphisms $h_i:\partial N(L_i)\rightarrow \partial N(K_i)$ such that $g(M)=g(E(L))-n$, where $M=E(L)\cup_{\cup h_i}(\bigcup E(K_i))$.\\

\begin{proof}

The case $n=1$ is just Proposition 6.4.  Using the method of \cite{schir}, we may construct a link $L=L_1\cup \cdots \cup L_n$ and a surface $S=S_1\cup S_2$ which is a disjoint union of two spheres in $S^3$, with the following properties:\\

\begin{enumerate}

\item $|S_j\cap L_i|=4$ for $j=1,2$, and $1\leq i \leq n$

\item For each component $B$ of $E(S)$, $E(L\cap B, B)$ is a handlebody.

\item $g(E(L))=4n-1$\\

\end{enumerate}

Let $K_i$ be a $2$-bridge knot with $2$-bridge sphere $S_i'$ in $S^3$, $1\leq i\leq n$, and let $h_i:\partial N(L_i)\rightarrow \partial N(K_i)$ be a homeomorphism which takes $S\cap \partial N(L_i)$ to $S_i'\cap \partial N(K_i)$.  Then for all $1\leq i \leq n$ the annuli $Fr(K_i\cap B')$ are primitive in the handlebodies $E(K_i\cap B',B')$ for each component $B'$ of $E(S_i')$, and so we can deduce that the surface $F=(S\cap E(L))\cup_{\cup h_i}(\bigcup (S_i'\cap E(K_i))$ is a Heegaard surface of $M$.  And it is easily verified that $g(F)=3n-1$.\\

\end{proof}

\end{proposition}

It is perhaps not surprising that the knots of Proposition 6.4 and the links of Proposition 6.5 were constructed (independently) by Nogueira and the author to serve as instances of knots and links which experience high levels of tunnel number degeneration.  In \cite{schir2}, the author has also shown that all of the knots and links constructed in \cite{schir}, which are there named {\em optimal} links, have the property that $g(E(L))\geq g(M)-n$, where $n$ is the number of components of $L$ and where $M$ is any manifold obtained via Dehn filling along an integral $n$-tuple of slopes of $\partial N(L)$.\\

\end{document}